\documentclass[a4paper,12pt]{amsart}
\usepackage{fullpage,verbatim}
\usepackage[english,francais]{babel}
\usepackage{epsf}
\usepackage[dvips]{epsfig}
\usepackage{amssymb}
\usepackage{amsmath}
\usepackage[T1]{fontenc}
\usepackage{wasysym}
\usepackage{aeguill}

\newtheorem{theorem}{Theorem}[section]

\newtheorem{lemma}[theorem]{Lemma}

\newcommand{\finpreuve}{\hspace{\stretch{1}}{$\square$}}

\newcommand{\card}{{\rm card}}
\newcommand{\R}{\mathbb{R}}
\newcommand{\E}{\mathbb{E}}
\renewcommand{\P}{\mathbb{P}}

\newcommand{\Z}{\mathbb{Z}}

\def\proof{\noindent{\bf Proof. }}

\def\remark{\noindent{\bf Remark. }}

\def\remarks{\noindent{\bf Remarks. }}

\def\proofof#1{\noindent{\bf Proof #1. }}
\renewcommand{\epsilon}{\varepsilon}
\def\btab{\begin{eqnarray*}}
\def\etab{\end{eqnarray*}}
\def\beq#1{\begin{equation}\label{#1}}
\def\eeq{\end{equation}}

\def\boulet{\noindent $\bullet$ }

\begin{document}

\selectlanguage{english}


\title[Continuous first-passage and greedy paths models]{Continuous first-passage percolation and continuous greedy paths model: linear growth}

{
\author{Jean-Baptiste Gou\'er\'e}
\author{R{\'e}gine Marchand}
\address{Laboratoire de Math{\'e}matiques, Applications et Physique
Math{\'e}matique d'Orl{\'e}ans UMR 6628\\ Universit{\'e} d'Orl{\'e}ans\\ B.P.
6759\\
 45067 Orl{\'e}ans Cedex 2 France}
\email{Jean-Baptiste.Gou\'er\'e@univ-orleans.fr}

\address{Institut Elie Cartan Nancy (math{\'e}matiques)\\
Universit{\'e} Henri Poincar{\'e} Nancy 1\\
Campus Scientifique, BP 239 \\
54506 Vandoeuvre-l{\`e}s-Nancy  Cedex France}
\email{Regine.Marchand@iecn.u-nancy.fr}
}

\def\motsclefs{First-passage percolation, greedy paths, boolean percolation, random growth.}


\begin{abstract}
We study a random growth model on $\R^d$ introduced by Deijfen.
This is a continuous first-passage percolation model. 
The growth occurs by means of spherical outbursts with random radii in the infected region. 
We aim at finding conditions on the distribution of the random radii to determine whether the growth of the process is linear or not. 
To do so, we compare this model with a continuous analogue of the greedy lattice paths model
and transpose results in the lattice setting to the continuous setting.
\end{abstract}

\maketitle 

\section{Introduction and statement of the main results}

We study a random growth model on $\R^d$, $d\ge 1$,  introduced by Deijfen in~\cite{Deijfen-modele}.
The model can be thought as describing the spread of an infection in a continuous medium. We fix an initially infected region $S_0$ in $\R^d$ (with positive Lebesgue measure) and a distribution $\mu$ on $(0, +\infty)$.
Let us denote by $S_t$ the random subset of $\R^d$ that corresponds to the infected region at time $t$ and by $|S_t|$ its Lebesgue measure. The random growth process $(S_t)_{t \ge 0}$ is a Markov process whose dynamics is as follows. 
Given $S_t$, we wait an exponentially distributed random time with mean $|S_t|^{-1}$.
We then add to $S_t$ a random ball, whose centre is chosen uniformly on $S_t$ and whose radius is chosen accordingly to the law~$\mu$. 

In~\cite{Deijfen-modele}, Deijfen proved an asymptotic shape result,
namely the almost sure convergence of $t^{-1}S_t$ toward a deterministic Euclidean ball.
This convergence holds as soon as the growth of $S_t$ is not superlinear.
She provided a sufficient condition for this behavior of the growth:
the boundedness of the support of $\mu$. 
This condition is weakened by Deijfen, H\"aggstr\"om and Bagley in \cite{Deijfen-Haggstrom-Bagley-competition} 
to the existence of an exponential moment for $\mu$.

In this work, we aim at finding a necessary and sufficient condition on the distribution $\mu$ of the radii of added balls 
to ensure the asymptotic shape result for this growth process. 
In Theorem~\ref{th-forme-asymptotique}, we prove that
\begin{itemize}
\item  If
\begin{equation} 
\label{intro-cs}
\int_0^{+\infty}\left(\int_x^{+\infty} r\mu(dr)\right)^{1/d}dx<+\infty,
\end{equation}
then the growth is not superlinear and therefore the asymptotic shape result holds.
\item On the other hand,  if
\begin{equation} 
\label{intro-cn}
\int_0^{+\infty} r^{d+1} \mu(dr)=+\infty,
\end{equation}
then the growth is superlinear and therefore the asymptotic shape result does not hold.
\end{itemize}
In dimension $d=1$, Conditions (\ref{intro-cs}) and (\ref{intro-cn}) are exclusive and thus give a necessary and sufficient condition 
for the linear growth. 
Unfortunately, in dimension $d \ge 2$, there is a gap between these two conditions that we did not manage to fill. 
Note however that if there exists $\epsilon>0$ such that
$$
\int_1^{+\infty} r^{d+1}(\ln r)^{d+\epsilon} \mu(dr)<+\infty
$$
then (\ref{intro-cs}) holds.
The gap is therefore reasonably sharp. 

\bigskip
To establish the sufficient condition  (\ref{intro-cs}),
we introduce and study a continuous analogue to the greedy lattice paths model introduced by
Cox, Gandolfi, Griffin and Kesten in~\cite{greedy-1}. 
In Theorem~\ref{th-greedy}, we give a necessary condition and a sufficient condition 
for the integrability of the supremum of mean weights of paths in the continuous setting.
Those results mimic similar ones in the discrete setting.

A comparison between Deijfen's model and those continuous greedy paths then enables us to conclude. 
Note that the gap between (\ref{intro-cs}) and (\ref{intro-cn}) comes directly from a similar gap for the continuous lattice paths model,
gap which is itself similar to the one existing for the greedy lattice paths model.

\bigskip
In the following, the dimension $d\ge 1$ is fixed. On $\R^d$, we denote by $\|.\|$ the Euclidean norm, 
and by $B_r$ the closed Euclidean ball centered at the origin with radius $r$.

\subsection{Deijfen's model}
\label{s-intro-deijfen}
Let us first recall the growth model introduced by Deijfen in~\cite{Deijfen-modele}.
Instead of using the original construction of the process, we use the construction, given later in~\cite{Deijfen-Haggstrom-coexistence}
by Deijfen and H\"aggstr\"om, which makes the analogy with first-passage percolation clearer.
We follow the presentation of Gou\'er\'e in~\cite{G-forme-territoire-competition}.

We fix a probability measure $\mu$ on $(0,+\infty)$.
We also fix $\chi$, a Poisson point process on $\R^d\times [0,+\infty) \times (0,+\infty)$ whose intensity
is the product of the Lebesgue measure on $\R^d\times [0,+\infty)$ by the distribution $\mu$ on $(0,+\infty)$. 

Let us consider the complete directed graph with vertex set $\R^d$.
We associate a passage time $\tau$ with each edge as follows:
\begin{enumerate}
\item For all $x\in\R^d$ we let $\tau(x,x)=0$.
\item For each point $(c,t,r)\in\chi$ -- where $c,t$ and $r$ respectively belong to $\R^d$, $[0,+\infty)$ and $(0,+\infty)$ --
and for each vertex $y\in (c+B_r)\setminus\{c\}$, we let $\tau(c,y)=t$.
\item For all edges $(x,y)$ to which we have not yet assigned any passage time, we let $\tau(x,y)=+\infty$.
\end{enumerate}
To say it in words, for a point $(c,t,r)$ in the Poisson process $\chi$, $t$ represents the time needed to travel from the center $c$ to each point of the ball centered in $c$ with radius $r$,  while outside balls, it takes an infinite time to travel. 

If $a$ and $b$ are two points of $\R^d$, we call path from $a$ to $b$ any finite sequence $\pi=(a=x_0,...,x_k=b)$ of distinct points of $\R^d$.
We denote by ${\mathcal C}(a,b)$ the set of such paths.
With each path $\pi=(x_0,...,x_k)$ we associate a passage time defined by:
$$
T(\pi)=\sum_{i=0}^{k-1} \tau(x_i,x_{i+1}).
$$
If $A$ is a subset of $\R^d$ and $x$ is a point of $\R^d$, we can then define the time $T(A,x)$
needed to cover $x$ starting from $A$, by:
$$
T(A,x)=\inf\{T(\pi) : a\in A, \; \pi \in {\mathcal C}(a,x)\}.
$$
Finally, if we start the process with the unit Euclidean ball $B$ centered at the origin, we can define the set $S_t$ of covered points at time $t$ by
$$
S_t = \{x \in \R^d : T(B,x)\le t\}.
$$

Relying on Kingman's subadditive ergodic theorem and using the isotropy of the model, one can establish the existence of a real $\lambda \ge 0$ such that
the following convergence holds almost surely:
$$
\lim_{\|x\|\to+\infty} \frac{T(B,x)}{\|x\|} = \lambda.
$$
This result is contained in the paper by Deijfen~\cite{Deijfen-modele}. 
The growth is linear if $\lambda$, which is the inverse of the speed, is positive.
In such a case, one can easily deduce the following asymptotic shape result: 
almost surely, for every $\epsilon>0$, for all large enough $t$, one has:
$$
B_{(1-\epsilon)/\lambda} \subset \frac{S_t}t \subset B_{(1+\epsilon)/\lambda}.
$$
As we need here this result under milder integrability assumptions that the ones used by Deijfen, 
we provide for completness a proof of this result in Appendix \ref{s-preuve-deijfen}.

To prove the positivity of $\lambda$ when $\mu$ admits an exponential moment, 
Bagley, Deijfen and H\"aggstr\"om introduce a new process that grows faster than Deijfen's one 
and whose non superlinear growth is easier to prove.
This new process can be roughly described as follows.
Assume that we have started with a set $A_0$ and that at time $t$ we have added balls $A_1,\dots, A_n$.
With each set $A_i$ we associate an exponential clock with mean $|A_i|^{-1}$.
These clocks are independent.
We wait for the first clock to ring.
If it is clock $i$, then we choose a point uniformly in $A_i$ and add a ball centered at this point with random radius. 
The projection on the first-coordinate axis of this new process is a one dimensional spatial branching process 
whose linear growth, when $\mu$ admits an exponential moment, is well-known.



Our conditions for the linear growth of Deijfen's model are  the following ones.

\begin{theorem} \label{th-forme-asymptotique}
1. If 
\begin{equation} \label{CS-vitesse}
\int_0^{+\infty}\left(\int_x^{+\infty} r\mu(dr)\right)^{1/d}dx<+\infty,
\end{equation}
then there exists $\lambda \in(0,+\infty)$ such that, almost surely, for all $\epsilon>0$, for all large enough $t$, one has:
$$
B_{(1-\epsilon)/\lambda} \subset \frac{S_t}t \subset B_{(1+\epsilon)/\lambda}.
$$
2. If
\begin{equation} \label{CN-vitesse}
\int_0^{+\infty} r^{d+1}\mu(dr)=+\infty,
\end{equation}
then, almost surely, for all $M>0$,  for all large enough $t$, one has:
$$
B_M \subset \frac{S_t}t.
$$
\end{theorem}

% ...

\remarks 
1. Note that in dimension $d=1$, if $\int_{(0,+\infty)}r^2\mu(dr)<+\infty$, one can compute explicitly the speed of the growth. 
With the notations of the theorem, one finds $\lambda=(\frac12 \int_0^{+\infty} r^2\mu (dr))^{-1}$.
We sketch a proof of this result in Appendix \ref{s-vitesse-dim1}.

2. 
To see how quantities scale,
let us prove that if one multiplies the radii by $2$, then one multiplies the speed by $2^{d+1}$.
Let us consider the following Poisson point processes:
\begin{eqnarray*}
\chi_1 & = & \{(c,t,2r), (c,t,r)\in\chi\}, \\
\chi_2 & = & \{(2^{-1}c,2^{-1}t,r), (c,t,r)\in\chi\}, \\
\chi_3 & = & \{(c,2^{-(d+1)}t,r), (c,t,r)\in\chi\}.
\end{eqnarray*}
With those points processes one can, as we have done with $\chi$, associate passage times to paths.
We denote them by $T_1,T_2$ and $T_3$.
Let $\pi=(x_0,\dots,x_n)$ be a path that originates in $0$. Then
$$
T_1(2\pi)\|2x_n\|^{-1}=T_2(\pi)\|x_n\|^{-1}, \quad T_2(\pi)\stackrel{\hbox{law}}{=}T_3(\pi), \quad T_3(\pi)=2^{-(d+1)}T(\pi).
$$
(The second property results of the equality in law of the Poisson point processes $\chi_2$ and $\chi_3$,
which itself results of the equality of their intensities.)
This gives the announced scaling.

\medskip
We would like now to explain how trying to rule out the possibility of an infinite speed in Deijfen's model 
naturally leads to introduce continuous greedy paths. We fix a small $\alpha>0$ and consider the "fast" balls:
\beq{def-xialpha}
\xi_\alpha=\{(c,r): \; \exists t \le \alpha r : \; (c,t,r) \in \xi\}.
\eeq
Roughly speaking, the speed in these balls is at least $1/\alpha$, while outside these "fast" balls, the speed is less than $1/\alpha$. 
More precisely, consider a path $\pi=(x_0,\dots,x_n)$ such that $T(\pi)$ is finite. 
Then, by the definition of $T$, every point $x_i$ in $\pi$ -- but the last -- 
 is the first coordinate of a point $(x_i,t_i,r_i)$ in $\chi$.
Moreover, we have $\|x_{i+1}-x_i\|\le r_i$ and
$$
T(\pi) =\sum_{i=0}^{n-1} t_i.
$$
If the path uses only "slow" balls, then
$$
T(\pi) =\sum_{i=0}^n t_i \ge \sum_{i=0}^n \alpha r_i \ge \sum_{i=0}^n \alpha \|x_{i+1}-x_i\| =\alpha |\pi|
$$
where $|\pi|$ is the length of the path $\pi$, that is the sum of the Euclidean length of its segments. 

If it also uses a "fast" ball $(x_i,t_i,r_i)$, the portion between $x_i$ and $x_{i+1}$ is traveled through at high speed. By considering that this 
portion is traveled through at infinite speed, we obtain
(the sums are over visited "fast" balls):
$$
T(\pi) \ge \alpha \left( |\pi|-\sum \|x_{i+1}-x_i\| \right) \ge \alpha \left( |\pi|-\sum r_i \right), \text{ or } 
\frac{T(\pi)}{|\pi|}\ge \alpha\left(1 -\frac{1}{|\pi|}\sum r_i \right).
$$
We are therefore led to bound from above the following kind of means:
$$
\frac{1}{|\pi|}\sum_{ x_i \in \pi, (x_i,r_i) \in \xi_\alpha} r_i.
$$
This motivates the introduction and the study of the continuous greedy paths discussed in Subsection~\ref{s-intro-greedy}. This crude link between 
Deijfen's model and continuous greedy paths will be precised in Subsection~\ref{s-lien}. The proof of Theorem~\ref{th-forme-asymptotique} 
is given in Subsection~\ref{s-th-forme-asymptotique}. 

\medskip
We end this subsection by pointing out a link with the Boolean model of continuum percolation and by giving some further intuition.
Fix $\alpha>0$ and denote by $\Sigma$ the union of the balls $c+B_r$, $(c,r)\in\xi_{\alpha}$
($\xi_{\alpha}$ is defined by (\ref{def-xialpha})).
This is the Boolean model of continuum percolation driven by $\xi_{\alpha}$.
We can define a first passage percolation process on the complete non-oriented graph of $\R^d$ as follows:
the time needed to travel along an edge $xy$ is $\alpha$ times the one-dimensional Lebesgue measure of $[x,y] \setminus \Sigma$.
In other words, one travels at speed $\alpha^{-1}$ outside $\Sigma$ and at speed $+\infty$ inside $\Sigma$.
By coupling, we can see that the speed in this first passage percolation process is larger than in Deijfen's model.

The intensity of $\xi_{\alpha}$ is the product of the Lebesgue measure on $\R^d$ by the measure $\nu_\alpha(dr)=\alpha r\mu(dr)$, which is also
$\alpha$ times the product of the Lebesgue measure on $\R^d$ by the measure $r\mu(dr)$.
By a result of Gou\'er\'e~\cite{G-percolation-continue}, one knows that the connected component of $\Sigma$ that contains the origin is almost
surely bounded for small enough $\alpha$ if and only if $\int r^d \nu_\alpha(dr)$ is finite, i.e.\ if and only if $\int r^{d+1} \mu(dr)$ is finite.

This suggests that, when $\int r^{d+1} \mu(dr)$ is finite, a constant and positive proportion of any path should lie outside $\Sigma$.
If this result were true, 
then the speed in the first passage percolation process associated with the Boolean model 
-- and therefore the speed in Deifen's model -- would be finite as soon as $\int r^{d+1} \mu(dr)$ is finite.
Unfortunately, we do not know whether the finiteness of $\int r^{d+1} \mu(dr)$ is sufficient or not 
to bound away from $0$ the proportion of length that any path spends outside $\Sigma$.

On the other hand, if $\int r^{d+1} \mu(dr)$ is infinite, "fast" balls with speed larger than $1/\alpha$ percolate for all $\alpha>0$.
This suggests that the speed in Deifjen's model is at least $1/\alpha$ for every $\alpha>0$,
which means that the growth is superlinear. 

\subsection{Greedy lattice paths model}
In this subsection, $d \ge 2$.

One first gives to points $c$ of $\Z^d$ i.i.d.\ positive random weights $r(c)$ with common law $\nu$.
A path is here a finite sequence of distinct points of $\Z^d$
such that the Euclidean distance between any two consecutive points is $1$ and the length of a path is naturally the sum of the Euclidean lengths of its segments.
With each path one associates a weight which is the sum of the weights of its points.
If $n$ is a positive integer, we denote by $A_n$ the supremum of the weights of all paths with length $n$ that originates in $0$.
In~\cite{greedy-1}, the authors show that if there exists a real $\epsilon>0$ such that:
$$
\int_1^{+\infty} r^d (\ln r)^{d+\epsilon} \nu(dr)<+\infty,
$$
then there exists a real $M<+\infty$ such that:
$$
\limsup_{n\to +\infty} \frac{A_n}{n} \le M \quad a.s.
$$
This result was improved in~\cite{greedy-2} by Gandolfi and Kesten: 
under the same condition, $A_n/n$ converges almost surely and in $L^1$ toward a finite constant. 
Martin, in~\cite{Martin-greedy}, obtains the same results under a weaker assumption and with a much simpler proof:\\
$$\text{ if } \int_0^{+\infty} \nu([r,+\infty))^{1/d}<+\infty \text{ then } A_n/n \text{ converges a.s. and in } L^1 \text{ to a finite constant}.$$

As an intermediate step, he shows that
\begin{equation} 
\label{resultat-lattice}
\text{if }\displaystyle \int_0^{\infty} \nu([r,+\infty))^{1/d}<+\infty \text{ then }\sup \E \left( \frac{A_n}n \right)<+\infty.
\end{equation}
Deriving this property in a continuous setting will turn out to be sufficient for our purpose.

On the other hand, from results in~\cite{greedy-1} and~\cite{Martin-greedy}, one knows that
$$\text{if } \int_0^{+\infty} r^d \nu(dr)=+\infty \text{ then }A_n/n \text{ almost surely goes to }+\infty.$$

\subsection{Continuous greedy paths}
\label{s-intro-greedy}
In our continuous analogue, the points of the lattice $\Z^d$ are replaced by the points of a homogeneous Poisson point process on $\R^d$. 
Fix a finite measure $\nu$ on $(0,+\infty)$, 
and consider a Poisson point process $\xi$ on $\R^d\times (0,+\infty)$
whose intensity is the product of the Lebesgue measure on $\R^d$ by the measure $\nu$. 
We denote by $\Xi$ the projection of $\xi$ on $\R^d$: 
the point process $\Xi$ is thus a Poisson point process on $\R^d$ with intensity $\nu((0,+\infty))$ times the Lebesgue measure on $\R^d$.
If $x$ is a point of $\Xi$ we denote by $r(x)$ the only positive real number such that $(x,r(x))$ belongs to $\xi$. 
Thus, given $\Xi$, the weights $(r(x))_{x \in \Xi}$ are i.i.d.\ with common law $\nu((0,+\infty))^{-1}\nu$ (if $\nu((0,+\infty))$ is positive). 
For points $x\in \R^d$ that are not in $\Xi$, we set $r(x)=0$.

A path is a finite sequence of distinct points of $\R^d$.
(In the lattice model, consecutive points of a path are required to be nearest neighbours.
We do not require such a condition in our continous model, which is therefore not an exact analogue of the lattice model.)
We denote by $|\pi|$ its length, that is the sum of the Euclidean length of its segments. 
We define the weight $A(\pi)$ of a path $\pi=(x_0,\dots,x_n)$ by:
$$
A(\pi)=\sum_{i=0}^n r(x_i).
$$
We are interested in the finiteness of the supremum $S$ of the mean weights of paths, defined as
\begin{equation}
\label{S}
S=\sup\left\{\frac{A(\pi)}{|\pi|}\right\},
\end{equation}
where the supremum is taken over all paths whose length is positive and that originates in~$0$.
In order to explicit the dependence of $S$ on $\xi$ we shall sometimes use the notation~$S(\xi)$.  
We also introduce, for $l>0$, 
\begin{equation}
\label{Sl}
S_l=\sup\left\{\frac{A(\pi)}{|\pi|}\right\},
\end{equation}
where the supremum is now taken over all paths whose length is larger than $l$ and that originates in~$0$.

As we will mainly consider paths starting from $0$, to avoid extra discussion on the status of the origin, 
we will always work on the full event $\{0 \not\in \Xi\}$.  
Let us notice that we do not change $S$ if we take the supremum over all paths 
$(x_0,\dots,x_n)$ such that, in addition to the previous requirements, $x_i$ belongs to $\Xi$ for all $i\ge 1$.
This can be seen by the triangular inequality.
We shall use this remark when convenient without further reference.


We state in the following theorem sufficient conditions for $S$ to be either integrable or a.s. infinite. 
These conditions are similar to the ones obtained for the discrete setting.
\begin{theorem} \label{th-greedy} Assume $d\ge 2$.\\
1. If $\displaystyle
\int_0^{+\infty} \nu([r,+\infty))^{1/d}dr<+\infty$ then $\mathbb{E}S<+\infty$. \\
2. If $\displaystyle
\int_0^{+\infty} r^d \nu(dr)=+\infty$
then $S$ is a.s.\ infinite.
\end{theorem}

\medskip
\remark If $d=1$, then $\E(S)$ is infinite as soon as $\nu((0,+\infty))$ is positive: 
the contribution of the first positive point already has an infinite mean. 
Indeed, denote by $X$ the smallest positive point of $\Xi$. 
This is an exponential random variable with parameter $\nu((0,+\infty))$.
Moreover $r(X)$ is distributed according to $\nu((0,+\infty))^{-1}\nu$ and is independent of $X$.
As a consequence, $\E(r(X)X^{-1})$ is infinite and therefore $\E(S)$ is infinite. 
We shall therefore be led to study $\E(\lim_{l \to +\infty}S_l)$ when $d=1$ in Subsection~\ref{s-proof-greedy1}.

\medskip
We conclude this subsection by giving some ideas of the proof.
The second item of the theorem is rather straightforward
(actually, the supremum of $r(x)\|x\|^{-1}$, $x\in\Xi$ is already a.s.\ infinite).
The proof of the first item follows the proof of the corresponding result
in the lattice setting by Martin~\cite{Martin-greedy}. 
The first step consists in studying the case where $\nu$ is the Dirac mass at point $1$,  in which case $S$ is integrable. 
For a general measure $\nu$, one then distinguishes between the contribution of the points $x\in \Xi$ according to the value of the radius $r(x)$. Using the fact that the supremum of a sum is less than the sum of the supremum, this allows to obtain the following upper bound:
$$\E S(\xi)  \le\int_0^{+\infty} \E S(\xi^{\lambda}) d\lambda, $$
where $\xi^{\lambda}=\{(c,1) : c\in \Xi \hbox{ and } r(c) \ge \lambda\}.
$
The point is then that, by a scaling argument for Poisson point processes, one can express $\E S(\xi^{\lambda})$ as the product of 
$\nu([\lambda, +\infty))^{1/d}$ by the expectancy of $S$ in the case where $\nu$ is the Dirac mass at point $1$. 
This scaling argument explains the role played by the dimension $d$. 
The proof is given in Subsection~\ref{s-proof-greedy}.

Note also that the link between 
the measure $\nu$ in the continuous greedy paths model and the distribution $\mu$ in Deijfen's model is presented in Subsection~\ref{s-lien} 
(and thus the link between the conditions in Theorem~\ref{th-greedy} and in Theorem~\ref{th-forme-asymptotique}).




\section{Proofs}

\subsection{Continuous greedy paths: Proof of Theorem~\ref{th-greedy}}
\label{s-proof-greedy}

We keep the notations and objects introduced in Subsection~\ref{s-intro-greedy}. 
We begin with the case of a deterministic radius equal to $1$ and we denote by $\delta_1$ the Dirac mass at point $1$.

\begin{lemma} 
\label{l-s-dirac-fini}
If $d\ge2$ and $\nu=\delta_1$, then $\mathbb{E}S<+\infty$.
\end{lemma}
\proof Let $\alpha_0>0$ be such that
$$
\forall \alpha \ge \alpha_0, \;\int_{\R^d} \exp(1-\alpha\|x\|)dx <1.
$$
Let $\alpha \ge \alpha_0$ be fixed and fix also an integer $k\ge 1$.
Let us denote by $B(k,\alpha)$ the set of all finite sequences $(x_1,\dots,x_k)$ of distinct points of $\Xi$ such that:
$$
k \ge \alpha \left( \|x_1\|+\sum_{i=2}^{k} \|x_i-x_{i-1}\|  \right).
$$
We have:
\begin{eqnarray*}
\mathbb{P}(B(k,\alpha)\neq \varnothing)
 & \le & \mathbb{E}(\card(B(k,\alpha))) \\
 & = & \int_{(\R^d)^k} \mathbf{1}_{k \ge \alpha \left( \|x_1\|+\sum_{i=2}^{k} \|x_i-x_{i-1}\|  \right)}dx_1\dots dx_k\\
 & \le & \int_{(\R^d)^k} \exp\left(k-\alpha\left( \|x_1\|+\sum_{i=2}^{k} \|x_i-x_{i-1}\|  \right)\right) dx_1\dots dx_k\\
 & = & \left( \int_{\R^d} \exp(1-\alpha\|x\|)dx\right)^k.
\end{eqnarray*}
Let us denote by $F(\alpha)$ the following event: 
there exists a path $\pi$ originating in $0$, with positive length and whose points, except $0$, belong to $\Xi$, 
such that the inequality $A(\pi) \ge \alpha |\pi|$ holds. 
Decomposing on the number of points in the path, we get
$$\{S>\alpha\} \subset F(\alpha)=\bigcup_{k=1}^{+\infty}\{B(k,\alpha)\neq\emptyset\},$$ and thus
\begin{eqnarray*}
\mathbb{P}(S>\alpha) & \le & \sum_{k=1}^{+\infty} \left(\int_{\R^d} \exp(1-\alpha\|x\|)dx\right)^k \\
& \le & \int_{\R^d} \exp(1-\alpha\|x\|)dx \left(1-\int_{\R^d} \exp(1-\alpha\|x\|)dx\right)^{-1} \\
& \le & \int_{\R^d} \exp(1-\alpha\|x\|)dx \left(1-\int_{\R^d} \exp(1-\alpha_0\|x\|)dx\right)^{-1} \\
& \le & \alpha^{-d}\int_{\R^d} \exp(1-\|x\|)dx \left(1-\int_{\R^d} \exp(1-\alpha_0\|x\|)dx\right)^{-1},
\end{eqnarray*}
which is an integrable function of $\alpha$ since $d\ge 2$. 
%
\finpreuve

\medskip
We next give the scaling argument announced in the introduction. 
Denote by $\xi_{\nu}$ a Poisson point process on $\R^d \times (0,+\infty)$ with intensity the product of Lebesgue's measure on $\R^d$ by the positive finite measure $\nu$ on $(0,+\infty)$.

\begin{lemma}
\label{echelle} For any $m>0$,
$S(\xi_{m\nu})$ has the same law as $m^{1/d}S(\xi_{\nu}).$
\end{lemma}

\proof
We just need to notice that the random set
$
\left\{\left(m^{1/d}c,r\right): \;  (c,r)\in  \xi_{m\nu}\right\}$
 is a Poisson point process on $\R^d \times (0,+\infty)$ with intensity the product of Lebesgue's measure on $\R^d$ by the positive finite measure $\nu$ on $(0,+\infty)$.
 \finpreuve

\bigskip
\proofof{of Theorem~\ref{th-greedy}}  \\
1. Assume that 
$$
\int_0^{+\infty} \nu([\lambda,+\infty))^{1/d}  d\lambda<+\infty.
$$ 
and let us prove that $\E S<+\infty$.

First, we need to make the process with the Dirac mass appear, by, in a certain manner, decomposing on the different values of the support of $\nu$. In fact, the useful way to do so is to use the classical trick $\displaystyle r=\int_0^{+\infty}\mathbf{1}_{r\ge \lambda}d\lambda$: we have
\begin{eqnarray}
 \mathbb{E}S
 & = & \mathbb{E}\left(\sup_{\pi=(0,x_1\dots,x_n)} \int_0^{+\infty} \frac{\sum_i \mathbf{1}_{r(x_i) \ge \lambda}}{|\pi|} d\lambda\right) \nonumber \\
 & \le & \int_0^{+\infty} \mathbb{E}\left(\sup_{\pi=(0,x_1,\dots,x_n)} \frac{\sum_i \mathbf{1}_{r(x_i) \ge \lambda}}{|\pi|} \right) d\lambda  \nonumber \\
 & \le & \int_0^{+\infty} \mathbb{E}S(\xi^{\lambda}) d\lambda,  \label{chaud}
\end{eqnarray}
where $\xi^{\lambda}$ is the point process on $\R^d\times(0,+\infty)$ defined by:
$$
\xi^{\lambda}=\{(c,1) : c\in \Xi \hbox{ and } r(c) \ge \lambda\}.
$$
Notice that $\xi^{\lambda}$ is a Poisson point process whose intensity is the product of the Lebesgue measure on $\R^d$
by $\nu([\lambda,+\infty))\delta_1$.

Then, we use the scaling property: if $\widetilde{\xi}$ is a Poisson point process on $\R^d\times(0,+\infty)$
whose intensity is the product of the Lebesgue measure on $\R^d$ by $\delta_1$, then
\begin{equation} \label{tuto}
\mathbb{E}S(\xi^{\lambda})=\nu([\lambda,+\infty))^{1/d} \mathbb{E}S(\widetilde{\xi}).
\end{equation}
Indeed, 
if $\nu([\lambda,+\infty))=0$ then the equality is straightforward, while if
 $\nu([\lambda,+\infty))\neq 0$, it is a simple application of Lemma~\ref{echelle}.

Finally, from (\ref{chaud}) and (\ref{tuto}) we get:
\begin{equation}
\label{ohlabellemajoration}
\mathbb{E} S\le \int_0^{+\infty} \nu([\lambda,+\infty))^{1/d} \mathbb{E}S(\widetilde{\xi}) d\lambda<+\infty
\end{equation}
by the integrability assumption on $\nu$ and  Lemma~\ref{l-s-dirac-fini}.


\bigskip
\noindent
2. Assume that 
$$\int_0^{+\infty}r^d \nu(dr)=+\infty$$
and let us prove that $S=+\infty$ a.s.

Let $M>0$ and  consider the following point process:
\begin{equation}
\label{troplourd}
\{(c,r) : (c,r)\in\xi, r \ge M\|c\| \hbox{ and } c\neq 0\}.
\end{equation}
The cardinal of this point process is distributed according to a Poisson law with parameter:
$$
\int_{\R^d} dc \int_0^{+\infty} 1_{r \ge M\|c\|} \nu(dr)=\int_0^{+\infty} |B_{rM^{-1}}| \nu(dr)=|B_{M^{-1}}|\int_0^{+\infty}r^d \nu(dr).
$$
By our assumption on $\nu$, this equals infinity and therefore this point process is almost surely non-empty.
But if $(c,r)$ is a point of this point process then, by considering the path $(0,c)$, we get $S\ge M$ a.s, which concludes the proof.
%
\finpreuve

\medskip

\subsection{Continuous greedy paths in dimension $d=1$}
\label{s-proof-greedy1}

To tackle the one dimensional case, where $\E (S)=+\infty$, we will rather use the asymptotic behaviour of $S_l$ stated in the following easy result. 

\begin{lemma} \label{Stilde-dim1} 
If $d=1$ and $\int_0^{+\infty} r\nu(dr)$ is finite, then:
$$
\lim_{l\to+\infty} S_l \le 2\int_0^{+\infty} r\nu(dr) \quad \hbox{ a.s.\ }
$$ 
\end{lemma}
\proof Let $l>0$ and $\pi=(x_0,\dots,x_n)\in \Pi_l$.
As $x_0=0$, all the $x_i$ belongs to $[-|\pi|,|\pi|]$.
We then have:
$$
A(\pi) \le \sum_{x\in\Xi \cap [-|\pi|,|\pi|]} r(x).
$$
Therefore:
\beq{brest}
\lim S_l \le \limsup \frac{1}{l} \sum_{x\in\Xi \cap [-l,l]} r(x).
\eeq
Recall the following:
\begin{itemize}
\item $\Xi$ is a Poisson point process with intensity $\nu((0,+\infty))$ times the Lebesgue measure.
\item Given $\Xi$, the sequence ${(r(x))}_{x\in \Xi}$ is an i.i.d.\ sequence of random variables distributed according to $\nu((0,+\infty))^{-1}\nu$.
\end{itemize}
Therefore, the right-hand side of (\ref{brest}) is a.s.\ 
$$
2\nu((0,+\infty)) \int_{(0,+\infty)} r\nu((0,+\infty))^{-1}\nu(dr).
$$
This concludes the proof. \finpreuve

%

\subsection{A link between Deijfen's Model and continuous greedy paths}
\label{s-lien}

Let us recall that the Poisson point process $\chi$,  driving Deijfen's model, has been introduced in Subsection~\ref{s-intro-deijfen}.
In this subsection, we fix a real $\alpha>0$.
We consider the following point process on $\R^d\times(0,+\infty)$:
$$
\xi_{\alpha}=\{(c,r):\exists t\le\alpha r : (c,t,r) \in \chi\}.
$$
In other words, $\xi_{\alpha}$ is the projection on $\R^d\times(0,+\infty)$ of the intersection of $\chi$
with the Borel set:
$$
\{(c,t,r)\in\R^d\times[0,+\infty)\times (0,+\infty) : t\le\alpha r\}.
$$
Let us notice that $\xi_{\alpha}$ is a Poisson point process on $\R^d\times(0,+\infty)$ whose intensity
is the product of the Lebesgue measure on $\R^d$ by the finite measure
$\nu_{\alpha}$ on $(0,+\infty)$ defined by:
$$
\nu_{\alpha}(dr)=\alpha r\mu(dr).
$$
We consider the continous greedy paths model driven by $\xi_{\alpha}$.
Except for the name of this point process, we keep the notations and objects defined in Subsections~\ref{s-intro-greedy} and~\ref{s-intro-deijfen}.

The set $\xi_{\alpha}$ corresponds to the "fast" balls, i.e.\ balls where the infection progresses with a speed larger than $1/\alpha$. 
The next lemma gives the link between the travel time in Deijfen's model 
and the functional $S$ in the continous greedy paths model driven by $\xi_{\alpha}$: 
roughly speaking, outside the balls in $\xi_{\alpha}$, 
the travel time between two points is at least $\alpha$ times the Euclidean distance between the two points, 
and the existence of "fast" balls gives a correction term controlled by $S(\xi_{\alpha})$.
\begin{lemma} 
\label{l-ppp-greedy}
For all vector $x\in\R^d\setminus B$, one has: 
$$
\frac{T(B,x)}{\|x\|} \ge \alpha\left(1-S_{\|x\|}(\xi_{\alpha})-\frac{1}{\|x\|}\right).
$$
\end{lemma}

\proof
Let $x\in\R^d\setminus B$.
Let $\pi=(x_0,\dots,x_n)$ be a path from $B$ to $x$.
(In other words, $x_0$ belongs to $B$ and $x_n$ equals $x$.)
As $\|x\|$ is positive, $S_{\|x\|} \le S$.
Therefore, in order to prove the lemma, it is sufficient to check the following inequality:
$$
\frac{T(\pi)}{\|x\|} \ge \alpha\left(1-S_{\|x\|}(\xi_{\alpha})-\frac{1}{\|x\|}\right).
$$
We assume that for every $i \in \{1,\dots,n\}$, $x_i \neq 0$
(otherwise, if $x_i=0$, one uses the inequality $T(\pi)\ge T(x_i,\dots,x_n)$ and works with the path $(x_i,\dots,x_n)$).

We extend $\pi$ in a path $\widetilde{\pi}$ starting from $0$ by adding if necessary a first point $x_{-1}=0$ to $\pi$.
As $x_n=x$, the length of the path $\widetilde{\pi}$ is at least $\|x\|$.
To establish the lemma, it is therefore sufficient to prove that:
$$
\frac{T(\pi)}{\|x\|} \ge \alpha\left(1-\frac{A(\widetilde{\pi})}{|\widetilde{\pi}|}-\frac{1}{\|x\|}\right).
$$

Let $i\in\{0,\dots,n-1\}$.
Let us show the following inequality:
\begin{equation} \label{e-lien}
\tau(x_i,x_{i+1}) \ge \alpha\|x_i-x_{i+1}\|-\alpha r(x_i).
\end{equation}
Remember that if there exists $r_i$ (a.s. unique) such that $(x_i,r_i) \in \xi_\alpha$ then $r(x_i)=r_i$, and $r(x_i)=0$ in any other case.
Three cases arise:
\begin{enumerate}
\item If $x_i=x_{i+1}$, then $ \tau(x_i,x_{i+1}) =0$ and $\alpha r(x_i) \ge 0$, thus  (\ref{e-lien}) holds.
\item If $\tau(x_i,x_{i+1})$ is infinite, (\ref{e-lien}) is obvious.
\item Otherwise, there exist $r_i,t_i$ such that $(x_i,t_i,r_i)$ belongs to $\chi$ and $\|x_i-x_{i+1}\| \le r_i$, which implies $\tau(x_i,x_{i+1}) = t_i$.
\begin{itemize}
\item If $t_i > \alpha r_i$, as $\alpha r(x_i)\ge 0$, (\ref{e-lien}) holds.
\item If $t_i \le \alpha r_i$, then $(x_i,r_i)\in \xi_\alpha$ and thus $r(x_i)=r_i$, which gives (\ref{e-lien}).
\end{itemize}
\end{enumerate}

As $\alpha r(x_n)$ is non-negative, summing  (\ref{e-lien}) for $i\in\{0,\dots,n-1\}$ implies that
$
T(\pi) \ge \alpha(|\pi|-A(\pi)).
$
As, moreover:
$$
\|x_0-x_{-1}\|-r(x_{-1})=\|x_0\|-r(0)\le 1,
$$
we obtain (whether $\widetilde{\pi}=(x_{-1},\dots,x_n)$ or $\widetilde{\pi}=(x_0,\dots,x_n)$) that
$
T(\pi) \ge \alpha(|\widetilde{\pi}|-A(\widetilde{\pi})-1).
$
From $|\widetilde{\pi}| \ge \|x_n-0\|=\|x\|$, we then deduce:
\begin{equation}
\label{olabel}
\frac{T(\pi)}{\|x\|} \ge \frac{T(\pi)}{|\widetilde{\pi}|}
\ge
\alpha\left(1-\frac{A(\widetilde{\pi})}{|\widetilde{\pi}|}-\frac{1}{|\widetilde{\pi}|}\right)
\ge
\alpha\left(1-\frac{A(\widetilde{\pi})}{|\widetilde{\pi}|}-\frac{1}{\|x\|}\right).
\end{equation}
The lemma follows. 
 
%
\finpreuve

\subsection{Deijfen's model: Proof of Theorem~\ref{th-forme-asymptotique}}
\label{s-th-forme-asymptotique}

Let us recall the following result from Deijfen~\cite{Deijfen-modele}.
As the result is not explicitely stated in~\cite{Deijfen-modele}, we provide a proof in Appendix \ref{s-preuve-deijfen}.

\begin{theorem} \label{th-vitesse}
There exists a constant $\lambda \ge 0$ such that the following convergence holds almost surely and in $L^1$:
$$
\lim_{\|x\|\to+\infty} \frac{T(B,x)}{\|x\|}=\lambda.
$$
If $\lambda>0$ then, almost surely, for all $\epsilon>0$ and for all large enough positive real $t$, one has:
$$
B_{\lambda^{-1}(1-\epsilon)t} \subset S_t \subset B_{\lambda^{-1}(1+\epsilon)t}.
$$
If $\lambda=0$ then, almost surely, for all $M>0$ and for all large enough positive real $t$, one has:
$$
B_{Mt} \subset S_t.
$$
\end{theorem}

The following lemma will enable us the prove the first part of Theorem~\ref{th-forme-asymptotique}.
We have stated it in such a way that its proof does not rely on Theorem~\ref{th-vitesse}.

\begin{lemma} \label{l-cs-fini}
$$\text{If }\int_0^{+\infty} \left(\int_x^{+\infty} r\mu(dr)\right)^{1/d} dx < +\infty, \text{ then } \liminf_{\|x\|\to+\infty} \frac{\E T(B,x)}{\|x\|} > 0.$$
\end{lemma}

\proof  Recall that measures $\nu_{\alpha}$, $\alpha>0$, are defined in Subsection~\ref{s-lien} and that:
\begin{equation} \label{totu}
\nu_{\alpha}(dr)=\alpha \nu_1(dr) = \alpha r\mu(dr).
\end{equation}
By assumption, the following condition holds:
$$
\int_0^{+\infty} \nu_1([x,+\infty))^{1/d}dx <+\infty.
$$

$\bullet$ Let us first consider the case $d\ge 2$.
By Theorem~\ref{th-greedy}, we then get that $\E S(\xi_1)$ is finite.
Using (\ref{totu}), Lemma~\ref{echelle} ensures that
$
\E S(\xi_{\alpha})=\alpha^{1/d} \E S(\xi_1).
$
Lemma~\ref{l-ppp-greedy} implies then that, for all real $\alpha>0$:
$$
\liminf_{\|x\|\to\infty} \frac{\E T(B,x)}{\|x\|} \ge \alpha(1-\E S(\xi_{\alpha}))=\alpha(1-\alpha^{1/d}\E S(\xi_{1})).
$$
But this quantity is positive as soon as $\alpha$ is small enough. 

$\bullet$ Let us assume now that $d=1$.
By Lemma~\ref{l-ppp-greedy} we get:
$$
\liminf_{\|x\|\to+\infty} \frac{T(B,x)}{\|x\|} \ge \alpha(1-\lim_{l \to +\infty} S_l(\xi_{\alpha})).
$$
By Lemma \ref{Stilde-dim1} we then get, almost surely:
$$
\liminf_{\|x\|\to+\infty} \frac{T(B,x)}{\|x\|} 
 \ge \alpha\left(1-2\int_0^{+\infty} r\nu_{\alpha}(dr)\right) 
 = \alpha\left(1-2\alpha \int_0^{+\infty} r^2\mu(dr)\right).
$$
By Fatou Lemma we get: 
$$
\liminf_{\|x\|\to+\infty} \frac{\E T(B,x)}{\|x\|} \ge \alpha\left(1-2\alpha \int_0^{+\infty} r^2\mu(dr)\right). 
$$
But this quantity is positive as soon as $\alpha$ is small enough. 
This ends the proof of the lemma. \finpreuve

\medskip
\remark By optimizing in $\alpha$, we get that, for $d\ge 2$,
$$\lambda \ge \frac{d^d}{(d+1)^{d+1}} \left( \frac1{\E S(\xi_{1})}\right)^d.$$
Using the bound in (\ref{ohlabellemajoration}), we obtain that, if $\widetilde{\xi}_d$ is a Poisson point process on $\R^d\times(0,+\infty)$
whose intensity is the product of the Lebesgue measure on $\R^d$ by $\delta_1$, then
 $$\E S(\xi_{1})\le  \mathbb{E}S(\widetilde{\xi}_d) \int_0^{+\infty} \left( \int_\lambda^{+\infty} r\mu(dr)\right)^{1/d}  d\lambda$$
 Finally, there exists a positive constant $C_d$ that only depends on the dimension $d$ such that:
$$
\lambda \ge C_d\left(\int_0^{+\infty} \left(\int_{\lambda}^{+\infty} r\mu(dr)\right)^{1/d}d\lambda\right)^{-d}.
$$
The result is still true for $d=1$.

\medskip
The following lemma will enable us the prove the second part of Theorem~\ref{th-forme-asymptotique}.

\begin{lemma}
 \label{l-cn-fini}
 $$\text{If }\int_0^{+\infty}r^{d+1} \mu(dr)=+\infty, \text{ then } \lim_{x\to\infty} \frac{T(B,x)}{\|x\|} = 0 \quad a.s \text{ and in } L^1.$$
\end{lemma}

\proof
By Theorem~\ref{th-vitesse},
one can fix a real $\beta>0$  such that, almost surely,  for all $x$ such that $\|x\|$ is large enough, the following inequality holds:
\begin{equation} \label{beta}
T(B,x) \le \beta \|x\|.
\end{equation}
Let $M > 0$ and write
$
A=\{(c,t,r) \in \R^d \times [0,+\infty) \times ]0,+\infty[ : t \le \|c\| \le rM^{-1}\}.
$
The cardinal of $\chi \cap A$ is distributed according to a Poisson law with parameter:
\begin{eqnarray*}
\lambda
 & = & \int_{\R^d} dc \int_{[0,+\infty)} dt \int_{]0,+\infty[} \mu(dr) 1_A(c,t,r) \\
 & = & \int_{]0,+\infty[} \mu(dr) \int_{B_{rM^{-1}}} dc \|c\|\\
 & = & \int_{]0,+\infty[} \mu(dr) r^{d+1} \int_{B_{M^{-1}}} dc \|c\|
  =  +\infty
\end{eqnarray*}
because of the assumption made about $\mu$.
Therefore, the cardinal of $\chi \cap A$ is almost surely infinite.
Let $s$ be a positive real and write
$
A_s = \{(c,t,r) \in A : \|c\| \le s\}.
$
The cardinal of $\chi \cap A_s$ is distributed according to a Poisson law with parameter:
$$
\lambda_s = \int_{\R^d} dc \int_{[0,+\infty)} dt \int_{]0,+\infty[} \nu(dr) 1_{A_s}(c,t,r) \le |B_s|s < +\infty.
$$
The cardinal of $\chi \cap A_s$ is therefore almost surely finite. The two preceding observations imply that, almost surely, there exists a sequence $(c_n,t_n,r_n)_n$
with values in $\chi \cap A$ such that $\|c_n\|$ goes to infinity.

As $Mc_n$ belongs to $c_n+B_{r_n}$, we have $\tau(c_n,Mc_n)=t_n$, which leads, using~(\ref{beta}) and the definition of $A$, to
$$T(B,Mc_n)\le T(B,c_n)+t_n\le \beta \|c_n\| + \|c_n\|=(1+\beta)\|c_n\|,$$
(for large enough $n$) and thus
$$
\frac{T(B,Mc_n)}{M\|c_n\|}\le \frac{1+\beta}{M},
\quad 
\text{ which implies } \quad 
\liminf_{\|x\|\to\infty} \frac{T(B,x)}{\|x\|} \le \frac{1+\beta}{M}.
$$
Theorem~\ref{th-vitesse} enables then to conclude the proof.
\finpreuve

\bigskip

\proofof{of Theorem~\ref{th-forme-asymptotique}} 
Just apply Theorem~\ref{th-vitesse}, Lemma~\ref{l-cs-fini} and Lemma~\ref{l-cn-fini}. 
\finpreuve

\appendix

\section{Computation of the speed in dimension $d=1$}
\label{s-vitesse-dim1}

We assume that $d=1$ and that $\int_0^{+\infty} r^2\mu(dr)$ is finite.
Our aim in this section is to sketch a proof of the following result: 
almost surely, for all $\epsilon>0$ and for all $t$ large enough, we have
$$
[-v(1-\epsilon)t,v(1-\epsilon)t] \subset S_t \subset [-v(1+\epsilon)t,v(1+\epsilon)t]
$$
where 
$$
v=\frac12\int_0^{+\infty} r^2\mu(dr)=\frac{\E(R^2)}{2}.
$$

We come back to the initial description of the process by Deijfen in \cite{Deijfen-modele}.
We keep the same point process $\chi$, but a point $(c,t,r) \in \chi$ has now the following interpretation: 
the ball centered at $c$ with radius $r$ becomes infected at time $t$ if its centre has been infected before time $t$. 
The equivalence of the two descriptions follows from the properties of the exponential law. 
In order to simplify the notation, instead of setting $S_0=[-1,1]$, we set $S_0=[-2,0]$.
We denote by ${(S_t)}_t$ the continuous first passage percolation process starting from $S_0$ and driven by $\chi$.
We denote by $\rho_t$ the rightmost point in $S_t$. 
By symmetry arguments, it is sufficient to prove that $\rho_t t^{-1}$ converges almost surely to $v$.
%

\boulet In this step, we show:
$$
\limsup_{t\to+\infty} \rho_t t^{-1} \le v, \hbox{ a.s.}
$$
In order to prove this result, we intoduce a process ${(\overline{S}_t)}_t$ that stochastically dominates ${(S_t)}_t$.
It has the same evolution by outbursts than ${(S_t)}_t$ - in particular it uses the same Poisson point process $\chi$ - 
but it starts from $\overline{S}_0=(-\infty,0]$. 
Denote by $\overline{\rho}_t$ the rightmost point in $\overline{S}_t$.
It suffices to show that $\overline{\rho}_t t^{-1}$ converges almost surely to $v$.

Denote by $\chi^+=\{(c,t,r) \in \chi: \; c\le 0, \; r+c>0\}$ the set of efficient balls. 
It is a Poisson point process with intensity measure:
$$
\mathbf{1}_{(-r,0]}(c)dc dt \mu(dr).
$$
Its projection on the $t$ coordinate is a homogeneous Poisson point process on $\R^+$ with finite intensity
$$
\int_0^{+\infty} \mu(dr) \int_{-r}^0 dc=\int_0^{+\infty}r\mu(dr)=\E R.
$$
Consider the first efficient outburst time $T_1$ (that is the leftmost point of the previous point process).
The random variable $T_1$ follows the exponential law with mean $(\E R)^{-1}$.
Let $(C_1,R_1)$ be such that $(C_1,T_1,R_1) \in \chi^+$. 
It is a random variable with law:
$$
(\E R)^{-1}\mathbf{1}_{(-r,0)}(c)dc \mu(dr).
$$
We have $\overline{\rho}_{T_1}=C_1+R_1$.
We can then compute the mean increase at time $T_1$:
$$
\E(\overline{\rho}_{T_1})=(\E R)^{-1}\int_0^{+\infty} \mu(dr) \int_{-r}^0 dc (c+r)=\frac{\E(R^2)}{2\E R}.
$$
Set $T_0=0$ and denote by ${(T_i)}_{i \ge 1}$ the sequence of efficient outburst times.
For example, 
$T_2$ is the $t$ coordinate of the point of $\{(c,t,r)\in\chi : c \le \overline{\rho}_{T_1}, \; r+c>\overline{\rho}_{T_1}, \; t\ge T_1\}$
whose $t$ coordinate is minimal.
The increments ${(\overline{\rho}_{T_i}-\overline{\rho}_{T_{i-1}})}_{i \ge 1}$ and ${(T_i-T_{i-1})}_{i \ge 1}$ 
are independent sequences of i.i.d.\ random variables.
Thus, 
$$
\frac{\overline{\rho}_{T_k}}k \to \E(\overline{\rho}_{T_1})=\frac{\E(R^2)}{2\E R} \quad \text{ and  } \quad   \frac{T_k}k \to \E(T_1)=\frac1{\E R}
$$
almost surely when $k$ goes to infinity.
As ${(T_i)}_{i \ge 1}$ goes almost surely to infinity, we can find for any $t>0$ a random $k$ such that $T_k \le t < T_{k+1}$.
For such a $k$, we have:
$$
\frac{\overline{\rho}_{T_k}}{T_{k+1}} \le \frac{\overline{\rho}_t}{t}\le \frac{ \overline{\rho}_{T_{k+1}}}{T_k}.
$$
This ends this step.

\boulet
In this step we show:
$$
\liminf_{t\to+\infty} \rho_t t^{-1} \ge v, \hbox{ a.s.}
$$
We introduce a process ${(\underline{S}_t)}_t$ that is stochastically dominated by ${(S_t)}_t$.
It also starts from $[-2,0]$, uses the same outburst process $\chi$, but increases only by the right side 
(the contribution of added balls on the left side are erased). 
Denote by $\underline{\rho}_t$ its rightmost point at time $t$. 
It suffices to prove:
$$
\liminf_{t\to+\infty} \underline{\rho}_t t^{-1} \ge v, \hbox{ a.s.}
$$

Set $T_0=0$.
We define $(C_1,T_1,R_1)$ as in the previous step.
The first efficient outburst for ${(\underline{S}_t)}_t$ occurs at time $T_1$ or later.
We have:
$$
\underline{\rho}_{T_1}=(C_1+R_1)\mathbf{1}_{C_1 \ge -2}.
$$
Let $T_2$ be the $t$-coordinate of the point of $\{(c,t,r)\in\chi : c \le \underline{\rho}_{T_1}, \; r+c>\underline{\rho}_{T_1}, \; t\ge T_1\}$
whose $t$ coordinate is minimal.
Let $C_2$ and $R_2$ be such that $(C_2,T_2,R_2)$ belongs to $\chi$.
We define in the same way $(C_i,T_i,R_i)$ for all $i\ge 3$.

Notice that ${(C_i-\underline{\rho}_{T_{i-1}},T_i-T_{i-1},R_i)}_{i \ge S1}$ is a sequence of i.i.d. random variables 
with the same law as $(C_1,T_1,R_1)$.
Moreover, for all $i\ge 1$, we have:
$$
\underline{\rho}_{T_i}-\underline{\rho}_{T_{i-1}}  = (C_i-\underline{\rho}_{T_{i-1}}+R_i)\mathbf{1}_{C_i \ge -2}.
$$
We write this as follows:
$$
\underline{\rho}_{T_i}-\underline{\rho}_{T_{i-1}} 
 =(C_i-\underline{\rho}_{T_{i-1}}+R_i)\mathbf{1}_{C_i-\underline{\rho}_{T_{i-1}} \ge -2-\underline{\rho}_{T_{i-1}}}.
$$
In particular, we have
$$
\underline{\rho}_{T_i}-\underline{\rho}_{T_{i-1}} 
 \ge (C_i-\underline{\rho}_{T_{i-1}}+R_i)\mathbf{1}_{C_i-\underline{\rho}_{T_{i-1}} \ge -2}
$$
from which we can conclude that $\underline{\rho}_{T_i}$ converges to $+\infty$.
Therefore, for all $s>0$, almost surely and for large enough $i$ we have:
$$
\underline{\rho}_{T_i}-\underline{\rho}_{T_{i-1}} 
 \ge (C_i-\underline{\rho}_{T_{i-1}}+R_i)\mathbf{1}_{C_i-\underline{\rho}_{T_{i-1}} \ge -s}.
$$
Consequently, for all $s>0$ and almost surely, we have:
\begin{eqnarray*}
\liminf_{i \to +\infty} \frac{\underline{\rho}_{T_i} }{i}
 &  \ge  & \liminf \frac{\sum_{j=1}^i (C_j-\underline{\rho}_{T_{j-1}}+R_j)\mathbf{1}_{C_j-\underline{\rho}_{T_{j-1}} \ge -s}}{i} \\
 & = & \E((C_1+R_1)\mathbf{1}_{C_1 \ge -s}).
\end{eqnarray*}
We let $s$ go to infinity and ends the proof as in the first step.

\boulet The result follows from the previous two steps. 
Let us notice that the proof of the previous step also allow us to prove again that the speed is infinite when $\E(R^2)$ is infinite.

\section{Proof of the asymptotic shape result (Theorem~\ref{th-vitesse})}
\label{s-preuve-deijfen}

Theorem~\ref{th-vitesse} is proved by Deijfen in \cite{Deijfen-modele}.
Nevertheless, we decided to provide a proof (which essentially follows Deijfen's one) for the following reasons:
\begin{itemize}
\item The theorem is not stated explicitely in \cite{Deijfen-modele}.
\item Deijfen works under stronger assumptions (but does not use them in the parts of the proof we are interested in).
\item We use the construction of the process given later in~\cite{Deijfen-Haggstrom-coexistence}. 
Our proof is therefore slightly simpler.  
\end{itemize}

The proof relies on Kingman's theorem.
Here is a statement of this theorem:

\begin{theorem} \label{th-kingman}
Suppose $(X_{m,n}, 0\le m<n)$ ($m$ and $n$ are integer) is a family of random variables satisfying:
\begin{enumerate}
\item For all integers $l,m,n$ such that $0\le l<m<n$, one has $X_{l,n}\le X_{l,m}+X_{m,n}$,
\item The distribution of $(X_{m+k,n+k}, 0\le m<n)$ does not depend on the integer $k$,
\item $\E(X_{0,1}^+)<\infty$ and there exists a real $c$ such that, for all natural integer $n$, one has
$\E(X_{0,n})\ge -cn$.
\end{enumerate}
Then
$$
\lim_{n\to\infty} \E(X_{0,n})n^{-1} \hbox{ exists and equals }\gamma=\inf_n \E(X_{0,n})n^{-1},
$$
$$
X:=\lim_{n\to\infty} X_{0,n}n^{-1} \hbox{ exists a.s. and in }L^1, \hbox{ and}
$$
$$
\E(X)=\gamma.
$$
If, for all $k\ge 1$, the stationary sequence $(X_{nk,(n+1)k}, n\ge 1)$ is ergodic, then $X=\gamma$ a.s.
\end{theorem}

We need a more general definition for passage times.
If $A$ and $C$ are two subsets of $\R^d$, we define $T(A,C)$,
the time needed to cover $C$ starting from $A$, by:
$${T}(A,C)=\sup_{c\in C} T(A,c).
$$
Let us begin with the following easy result:
\begin{lemma} \label{it-DHB} 
If $A, C$ and $D$ are mesurable subsets of $\R^d$, then:
$$
{T}(A,D)\le {T}(A,C)+{T}(C,D).
$$
\end{lemma}
\proof Assume that the right hand side of the inequation stated in the lemma is finite
(otherwise the result is obvious).
Let $d\in D$ and $\epsilon>0$.
Fix $c\in C$ and $\pi_2\in {\mathcal C}(c,d)$ be such that ${T}(r_2)\le{T}(C,D)+\epsilon$.
Now, fix $a\in A$ and $r_1\in {\mathcal C}(a,c)$ such that ${T}(r_1)\le{T}(A,C)+\epsilon$.
If we concatenate $r_1$ and $r_2$, we get an element $r\in {\mathcal C}(a,d)$ such that
${T}(r)\le {T}(C,D)+{T}(A,C)+2\epsilon$.
The lemma follows. \finpreuve

In the following lemma, we give bounds on quantities such as $T(B,B+x)$.
The proof essentially relies on the observation that such quantities can be bounded above by the sum of $\lceil \|x\| \rceil$
independent exponential random variables. For $\gamma>0$, we denote by $B_\gamma$ the euclidean ball centered in the origin with radius $\gamma$.

\begin{lemma} \label{GD0}
Let $\gamma>0$.
There exists $C,D>0$ such that the following inequality holds for all $x\in\R^d$:
$$
\E(T(B_{\gamma},x+B_{\gamma}))\le C\|x\|+D.
$$
Moreover there exist $a,b,c>0$ such that the following inequality holds for all $r>0$:
$$
\P(T(B_{\gamma}, B_r)\ge ar)\le b\exp(-cr).
$$
\end{lemma}
\proof 

\boulet Let $\gamma_0>0$ be such that $\nu([5\gamma_0,+\infty))$ is positive.
Let us notice that it is sufficient to prove the lemma for $\gamma\le \gamma_0$.
This is a consequence of the following two observations:
\begin{itemize}
\item[-] Let $\gamma \ge \gamma_0$. 
If $A$ be a finite subset of $\R^d$ such that $B_{\gamma}$ is contained in $A+B_{\gamma_0}$ then,
for all $x\in\R^d$, we have:
\btab
T(B_{\gamma},x+B_{\gamma}) 
 & \le & T\left(B_{\gamma_0},\bigcup_{y\in A} x+y+B_{\gamma_0}\right) \\
 & = & \max_{y\in A} T(B_{\gamma_0},x+y+B_{\gamma_0}) \\
 & \le & \sum_{y\in A} T(B_{\gamma_0},x+y+B_{\gamma_0}).
\etab
\item[-] For all $r>0$, the map $\gamma \mapsto T(B_{\gamma}, B_r)$ is non-increasing.
\end{itemize}

\boulet Hencesoforth, we assume that $\gamma \in (0, \gamma_0]$.
Set 
$$\lambda=|B_{\gamma}| \, \nu([5\gamma,+\infty))>0.$$
Let ${(T_n)}_{n\ge 0}$ be a sequence of i.i.d.\ exponential random variables with parameter $\lambda$. Let us prove that, for all $x\in\R^d$
\begin{equation}
\label{dom}
T(B_{\gamma},x+B_{\gamma}) \text{ is stochastically dominated by }T_0+\dots+T_n, \text{ with } n=\left\lfloor\frac{\|x\|}{3\gamma}\right\rfloor.
\end{equation}
Let $x\in\R^d \setminus \{0\}$ and let $n$ be defined in (\ref{dom}).
For all $i\in\{0,\dots,n\}$ we write
$$
x_i= \frac{3\gamma i}{\|x\|} x
$$ 
and define the set 
$
C_i=(x_i+B_{\gamma}) \times [0,+\infty) \times [5\gamma, +\infty).
$
For all $i$, we denote by $U_i$ the time coordinate of the point of $\chi \cap C_i$ whose time coordinate is minimal.
The $U_i$ have the same distribution, namely the exponential distribution with parameter $\lambda$ 
Moreover, as the sets $C_i$ are disjoint, the random variables $U_i$ are independent.
Write $x_{n+1}=x$.
Notice that, for all $c \in x_i+B_{\gamma}$ and all $r\ge 5\gamma$, the set $x_{i+1}+B_{\gamma}$ is contained in the set $c+B_r$.
Therefore, we have:
$$
T(x_i+B_{\gamma},x_{i+1}+B_{\gamma}) \le U_i
$$
and then, by the triangular inequality (Lemma \ref{it-DHB}):
$$
T(B_{\gamma},x+B_{\gamma}) \le U_0+\dots+U_n,
$$
which proves (\ref{dom}). The first part of the lemma follows immediately.

\boulet In this step, we prove the existence of $a_1,b_1,s_1>0$ such that the following inequality holds for all $x\in\R^d$
and all real $s \ge s_1$:
\beq{tour}
\P\big(T(B_{\gamma}, x+B_{\gamma}) \ge s\|x\|\big)\le a_1\exp(-b_1s\|x\|).
\eeq
Let $x\in\R^d$ and $s>0$.
By (\ref{dom}), 
$
\P\big(T(B_{\gamma},x+B_{\gamma}) \ge s\|x\|\big) \le \P\big(T_0+\dots+T_n \ge s\|x\|\big).  
$
Therefore 
$$
\P(T(B_{\gamma},x+B_{\gamma}) \ge s\|x\|) \le a_1^{n+1} \exp\left(-\frac{\lambda s\|x\|}{2}\right), \text{ where }a_1=\E\left(\exp\left(\frac{\lambda T_0}{2}\right)\right)<+\infty.
$$
Let $s_1>0$ be such that, for all $s\ge s_1$, the following inequality holds:
$$
a_1\exp\left(-\frac{3 \lambda s\gamma}{4}\right)<1.
$$
For such $s$ we then have:
\btab
\P\big(T(B_{\gamma},x+B_{\gamma}) \ge s\|x\|\big)
 & \le & a_1^{n+1} \exp\left(-\frac{\lambda s \|x\| }{4}\right) \exp\left(-\frac{\lambda s\|x\|}{4}\right) \\ 
 & \le & a_1^{n+1} \exp\left(-\frac{3 \lambda s\gamma n }{4}\right) \exp\left(-\frac{\lambda s\|x\|}{4}\right) \\ 
 & \le & a_1 \exp\left(-\frac{\lambda s\|x\|}{4}\right).
\etab
Choosing $b_1=\lambda/4$ we see that (\ref{tour}) holds.

\boulet
Let $M>0$ be such that, for all $r>0$, there exists a set $\Sigma \subset B_r \setminus\{0\}$ such that:
\beq{debussy}
\card(\Sigma) \le Mr^d \hbox{ and } B_r \subset \Sigma+B_{\gamma}.
\eeq
Let $r>0$. 
Let $\Sigma \subset B_r$ be such that (\ref{debussy}) holds.
We then have:
\btab
\P(T(B_{\gamma},B_r) \ge s_1r) 
 & \le & \P\left(\max_{x\in\Sigma} T(B_{\gamma},x+B_{\gamma}) \ge s_1r \right) \\
 & \le & \sum_{x\in\Sigma} \P(T(B_{\gamma},x+B_{\gamma}) \ge s_1r ) \\
 & \le & \sum_{x\in\Sigma} \P\left(T(B_{\gamma},x+B_{\gamma}) \ge \frac{s_1r}{\|x\|}\|x\| \right).
\etab
As $\Sigma$ is contained in $B_r\setminus\{0\}$, we have $s_1r\|x\|^{-1} \ge s_1$ for all $x\in\Sigma$.
By (\ref{tour}), we then get:
$$
\P(T(B_{\gamma},B_r) \ge s_1r) 
  \le  \sum_{x\in\Sigma} a_1\exp(-b_1 s_1r) \\
  \le Mr^da_1\exp(-b_1 s_1r) ,
$$
which proves the second point of the lemma. \finpreuve

\bigskip
From the previous lemma and Borel-Cantelli's lemma, we deduce the following result.
\begin{lemma} \label{GD}
Let $\beta,\epsilon>0$. There exists $\alpha>0$ such that, on a full probability event, the following holds:
$$
\hbox{The set } \Big{\{}x \in \R^d : T(x+B_{\beta},x+B_{\alpha \|x\|} ) \ge \epsilon \|x\|\Big{\}} \hbox{ is bounded}.
$$
\end{lemma}
\proof Let $\alpha, \beta, \epsilon>0$.

\boulet Let $\gamma>0$ be such that any ball with radius $\beta$ contains a ball with radius $\gamma$ centered at a point of $\gamma \Z^d$.
Let us prove that 
\begin{eqnarray}
\label{bounded}
&& \text{if the set }H=\Big{\{}y \in \gamma\Z^d : T(y+B_{\gamma},y+B_{3\alpha\|y\|}) \ge \epsilon \|y\|/2\Big{\}} \text{ is finite,}
\nonumber \\
&& \text{then the set }  G=\Big{\{}x \in \R^d : T(x+B_{\beta},x+B_{\alpha \|x\|}) \ge \epsilon \|x\|\Big{\}}\text{ is bounded.} 
\end{eqnarray}
Let $x\in\R^d$ be such that:
\beq{eq-con}
\|x\|\ge\max(3\beta,\beta+\beta\alpha^{-1}).
\eeq
Let $y\in \gamma\Z^d$ be such that $y+B_{\gamma}$ is contained in $x+B_{\beta}$.
We have, for all $z\in x+B_{\alpha\|x\|}$:
$$
\|z-y\|\le \|z-x\|+\|x-y\| \le \alpha\|x\|+\beta \le  \alpha(\|y\|+\beta+\beta\alpha^{-1}).
$$
By (\ref{eq-con}), we have $\|y\|\ge \|x\|-\beta \ge \max(2\beta,\beta\alpha^{-1})$, and thus
$
\|z-y\|\le 3\alpha \|y\|.
$
In other words, we have:
$$
x+B_{\alpha\|x\|}\subset y+B_{3\alpha\|y\|}.
$$
Therefore, if $T(x+B_{\beta},x+B_{\alpha \|x\|} +x)$ is greater that $\epsilon \|x\|$, 
then $T(y+B_{\gamma},y+B_{3\alpha\|y\|})$ is also greater that $\epsilon \|x\|$.
Moreover, as $\|y\| \ge 2\beta$, we have
$
\|x\| \ge \|y\|-\beta \ge \|y\|/2,
$
therefore $T(y+B_{\gamma},y+B_{3\alpha\|y\|})$ is greater that than $\epsilon \|y\|/2$. 
This proves (\ref{bounded}).

\boulet We now prove that, for a small enough $\alpha>0$, the set $H$ is almost surely finite:
This will conclude the proof of the lemma.
Let $a,b,c>0$ be the real numbers given by the second part of Lemma \ref{GD0}.
Set $\alpha=\epsilon/(6a)$.
We have:
\btab
\E(\card(H))
 & = & \sum_{y\in\gamma\Z^d} \P\big(T(y+B_{\gamma},y+B_{3\alpha\|y\|}) \ge \epsilon \|y\|/2 \big) \\
 & = & \sum_{y\in\gamma\Z^d} \P\big(T(B_{\gamma},B_{3\alpha\|y\|}) \ge \epsilon \|y\|/2 \big) \\
 & = & \sum_{y\in\gamma\Z^d} \P\big(T(B_{\gamma},B_{3\alpha\|y\|}) \ge a(3\alpha\|y\|) \big) \\
 & \le & \sum_{y\in\gamma\Z^d}b\exp(-3c\alpha\|y\|) <+\infty,
\etab
which concludes the proof. \finpreuve

%

\bigskip
The following result is a consequence of Kingman's theorem and of the previous lemmas.
If we consider the behaviour of $T(B,\cdot+B)$ and not directly the behaviour of $T(B,\cdot)$,
it is because $T(B,\cdot+B)$ has the following subbadditive property: 
$$
T(B,x+y+B) \le T(B,x+B)+T(x+B,x+y+B).
$$

\begin{lemma} \label{step1}
There exists $\lambda \ge 0$ such that the following convergence holds almost surely and in $L^1$:
$$
\lim_{\|x\|\to+\infty} \frac{T(B,x+B)}{\|x\|} = \lambda.
$$
\end{lemma}
\proof
\boulet Let us denote by $S=\{x\in\R^d:\|x\|=1\}$ the unit Euclidean sphere.
We first prove the existence of a real $\lambda \ge 0$ such that,
for all vector $x\in S$, the following convergence holds almost surely and in $L^1$:
$$
\lim_{n\to+\infty} \frac{T(B,nx+B)}{n} = \lambda.
$$
Let $x\in S$.
For all integers $m,n$ we set $X_{m,n}=T(B+mx,B+nx)$.
The first condition of Kingman's theorem is satisfied thanks to Lemma~\ref{it-DHB}.
The second one and the ergodicity one are consequences of related properties of the Poisson point process $\chi$.
The third one is satisfied because of the first part of Lemma \ref{GD0}.
We therefore have the convergence of $T(B,kx+B)k^{-1}$ (almost surely and in $L^1$) toward a non-negative deterministic real.
Because of the isotropy property of the model (which is a consequence of the related property of $\chi$), this limit does not depend on $x\in S$.
We denote it by $\lambda$.

%

\boulet 
Let $\epsilon>0$.
Let $C$ and $D$ be the real numbers given by Lemma \ref{GD0} with $\gamma=1$.
Let $\alpha$ be the positive real given by Lemma \ref{GD} with $\beta=1$.
We can assume that $\alpha$ is also such that the inequality $C\alpha\le \epsilon$ holds.

Let $\Sigma$ be a finite subset of $S$ such that $\Sigma+B_{\alpha/4}$ contains $S$.
Let $y\in \R^d$ be such that $\|y\|$ is large ($\|y\| \ge \max(1,8\alpha^{-1})$ is enough).
Set $n=\lfloor \|y\|\rfloor$. 
Let $x\in \Sigma$ be such that $y\|y\|^{-1}$ belongs to $x+B_{\alpha/4}$.
We can write: 
$$
\frac{T(B,y+B)}{\|y\|}-\lambda=I+J+K
$$
where
$$
I=\frac{T(B,y+B)}{\|y\|}-\frac{T(B,nx+B)}{\|y\|},
J=\frac{n}{\|y\|}\left(\frac{T(B,nx+B)}{n}-\lambda\right)
\text{ and }
K=\lambda\left(\frac{n}{\|y\|}-1\right).
$$
Using the triangular inequality satisfied by $T$ (Lemma \ref{it-DHB}) and the inequality $\|nx\|\le \|y\|$ (recall that $\|x\|=1$), we get:
\btab
|I| 
 & \le &  \max\left(\frac{T(nx+B,y+B)}{\|y\|}, \frac{T(y+B,nx+B)}{\|y\|}\right) \\
 & \le &  \max\left(\frac{T(nx+B,y+B)}{\|nx\|}, \frac{T(y+B,nx+B)}{\|y\|}\right) 
\etab
Notice that:
\beq{pppp}
\|nx-y\|+1 =\Big{\|}(nx-\|y\|x)+\|y\|(x-y\|y\|^{-1})\Big{\|} \le 2 + \alpha\|y\|/4 \le \alpha\|y\|/2
\eeq
(recall that we have $\|y\|\ge 8\alpha^{-1}$). 
Notice also the inequality $\|nx\|\ge \|y\|-1$ and therefore:
$$
\|nx-y\|+1 \le \alpha\|nx\|/2+\alpha/2 \le \alpha\|nx\|
$$
(recall that we have $\|y\|\ge 1$ and then $\|nx\|\ge 1$).
We then have:
$$
|I| \le \max\left(\frac{T(nx+B,nx+B_{\alpha \|nx\|})}{\|nx\|}, \frac{T(y+B,y+B_{\alpha \|y\|})}{\|y\|}\right).
$$
By definition of $\alpha$, we then get that, almost surely, if $y$ is large enough (which implies that $\|nx\|$ is also large), 
the inequality $|I|\le \epsilon$ holds.
Moreover, we have:
$$
\E(|I|) \le \E\left(\frac{T(nx+B,y+B)}{\|y\|}\right)+\E\left(\frac{T(y+B,nx+B)}{\|y\|}\right).
$$
By stationarity and by definition of $C$ and $D$ and by (\ref{pppp}) we get:
$$
\E(|I|) \le \frac{2C\|y-nx\|+2D}{\|y\|} \le C\alpha+\frac{2D}{\|y\|}.    
$$
As we have assumed $C\alpha \le \epsilon$, we get that $\E(|I|)$ is less than $2\epsilon$ if $\|y\|$ is large enough.

We have:
$$
|J| \le \left|\frac{T(B,nx+B)}{n}-\lambda\right|.
$$
Recall that $n\ge \|y\|-1$ and that $x$ belongs to the finite set $\Sigma$. 
By the first step we then get that almost surely, if $\|y\|$ is large enough, then $|J|\le \epsilon$.
We also get that if $\|y\|$ is large enough, then $\E(|J|)\le\epsilon$.
As moreover $|K|\le \lambda\|y\|^{-1}$, 
this concludes the proof. \finpreuve

\bigskip In the following Lemma we assume that the radius of the added balls is bounded from below by $\beta>0$.
As a consequence, when $S_t$ reaches a point $x$, $S_t$ covers a whole ball with radius $\beta$ that contains $x$.
From this observation and Lemma \ref{GD} it is easy to see that $T(B,x)$ and $T(B,x+B)$ have the same asymptotic behaviour. 

\begin{lemma} \label{step2}
Assume the existence of $\beta>0$ such that $\mu([\beta,+\infty))=1$.
Then there exists a real $\lambda \ge 0$ such that the following almost sure convergences hold:
$$
\lim_{\|x\|\to+\infty} \frac{T(B,x+B)}{\|x\|}=\lim_{\|x\|\to+\infty} \frac{T(B,x)}{\|x\|} = \lambda.
$$
\end{lemma}
\proof
Let $\epsilon>0$.
Let $\alpha$ be the real number given by Lemma \ref{GD}.
Let $M$ be an almost surely finite random variable such that the set:
$$
\Big{\{}x \in \R^d : T(B_{\beta}+x,B_{\alpha \|x\|} +x) \ge \epsilon \|x\|\Big{\}}
$$
is included in $B_M$.
Let $\lambda$ be the real number given by Lemma \ref{step1}.
Let $x\in\R^d$ be such that $\|x\|$ is large
($\|x\|> \max(1,M+\beta,\alpha^{-1}+\beta\alpha^{-1}+\beta)$ is enough).
Let us prove the following inequality:
\beq{lessive}
T(B,x+B) \le T(B,x)+1+\epsilon\beta+\epsilon \|x\|.
\eeq
We may assume that $T(B,x)$ is finite.
By definition of $T(B,x)$, there exists a path $\pi=(c_0,\dots,c_n)$ from a point of $B$ to $x$ such that $T(\pi)$ is less than $T(B,x)+1$.
As $\|x\|> 1$ we have $n\ge 1$.
As $T(B,x)$ is finite, $T(\pi)$ is finite and therefore $\tau(c_{n-1},x)$ is finite.
As a consequence, $x$ belongs to $B_{r(c_{n-1})}+c_{n-1}$.
But we have $r(c_{n-1})\ge \beta$.
Therefore there exists $y\in\R^d$ such that the following holds: $\|y-x\|\le\beta$ and $y+B_{\beta} \subset c_{n-1}+B_{r(c_{n-1})}$.
But for all $z\in c_{n-1}+B_{r(c_{n-1})}$ we can consider the path $\pi(z)=(c_0,\dots,c_{n-1},z)$.
We then have $T(\pi(z))=T(\pi)$ and then: 
$$
T(B,c_{n-1}+B_{r(c_{n-1})}) \le T(\pi) \le T(B,x)+1
$$
and therefore
$$
T(B,y+B_{\beta}) \le T(B,x)+1.
$$
From $\|x\| > M+\beta$ and $\|y-x\|\le\beta$ we deduce $\|y\| > M$.
By definition of $M$, we then have:
$$
T(y+B_{\beta},y+B_{\alpha \|y\|} ) \le \epsilon \|y\|.
$$
From the previous inequalities and the triangular inequality (see Lemma \ref{it-DHB}), we get:
$$
T(B,y+B_{\alpha \|y\|} ) \le T(B,x)+1+\epsilon \|y\| \le T(B,x)+1+\epsilon\beta+\epsilon\|x\|.
$$
From $\|x\| \ge \alpha^{-1}+\beta\alpha^{-1}+\beta$ and $\|y-x\|\le\beta$ we get that the ball $x+B$ is contained in the ball $y+B_{\alpha \|y\|}$.
Therefore (\ref{lessive}) holds.

We then have, almost surely:
$$
\liminf \frac{T(B,x+B)}{\|x\|} \le \liminf \frac{T(B,x)}{\|x\|}+\epsilon.
$$ 
Therefore, almost surely:
$$
\liminf \frac{T(B,x+B)}{\|x\|}-\epsilon \le \liminf \frac{T(B,x)}{\|x\|} \le \limsup \frac{T(B,x)}{\|x\|} \le \limsup \frac{T(B,x+B)}{\|x\|}
$$
and then, almost surely:
$$
\lambda -\epsilon \le \liminf \frac{T(B,x)}{\|x\|} \le \limsup \frac{T(B,x)}{\|x\|} \le \lambda.
$$
The lemma follows.
\finpreuve

\bigskip In the following Lemma, we get the asymptotic behaviour of $T(B,x)$ in the general case.
The idea is to compare the passage times $T(B,x)$ associated with different laws for the radius of the added balls, some of which fulfilling
the assumptions of the previous lemma.

\begin{lemma} \label{step3}
There exists a real $\lambda \ge 0$ such that the following almost sure convergences hold:
$$
\lim_{\|x\|\to+\infty} \frac{T(B,B+x)}{\|x\|}=\lim_{\|x\|\to+\infty} \frac{T(B,x)}{\|x\|} = \lambda.
$$
\end{lemma}
\proof
Let $\beta>0$ be small enough to ensure $\mu([\beta,+\infty))>0$.

\boulet Let us consider the following point process:
$$
\chi_-=\{(c,t,r) : (c,t,r) \in \chi \hbox{ and } r \ge \beta\}. 
$$
Let us denote by $T_-(\cdot,\cdot)$ the passage times associated with $\chi_-$ 
(in the same way as the $T(\cdot,\cdot)$ are associated with $\chi$).
We have, for all $x\in\R^d$ :
\beq{couche1}
T(B,x+B) \le T_-(B,x+B).
\eeq

\boulet Let us now consider the following point process:
$$
\widetilde{\chi_+} = \{(c,t \mu([\beta,+\infty),r) : (c,t,r) \in  \chi_-\}. 
$$
Its intensity is the product of the Lebesgue measure on $\R^d\times [0,+\infty)$ 
by the distribution $\mu_+$ on $(0,+\infty)$ defined by:
$$
\mu_+(A)=\frac{\mu(A \cap [\beta,+\infty))}{\mu([\beta,+\infty))}. 
$$
As the distribution $\mu_+$ stochastically dominates the distribution $\mu$, 
there exists a point process $\chi_+$ such that:
\begin{itemize}
\item $\chi_+\stackrel{\hbox{law}}{=}\widetilde{\chi_+}$.
\item For all $(c,t,r)$ in $\chi$, there exists $r' \ge r$ such that $(c,t,r')$ belongs to $\chi_+$.
\end{itemize}
We denote by $T_+$ and $\widetilde{T_+}$ the associated passage times.
We have, for all $x\in\R^d$ :
\beq{couche2}
T_+(B,x) \le T(B,x) \le T(B,x+B).
\eeq

\boulet Let us denote by $\lambda(\beta)$ the real given by Lemma \ref{step2} for the passage times $T_+(\cdot,\cdot)$.
We have the following almost sure convergences:
$$
\lim_{\|x\|\to+\infty} \frac{T_+(B,x)}{\|x\|} = \lambda(\beta)
\text{ and }
\lim_{\|x\|\to+\infty} \frac{T_+(B,x+B)}{\|x\|} = \lambda(\beta).
$$
As $\chi_+$ and $\widetilde{\chi_+}$ have the same laws the previous convergences still hold if we replace $T_+$ by $\widetilde{T_+}$.
Using the definition of $\widetilde{\chi_+}$, we then get the following almost sure convergences:
$$
\lim_{\|x\|\to+\infty} \frac{T_-(B,x+B)}{\|x\|} = \frac{\lambda(\beta)}{\mu([\beta,+\infty))}\ge \lambda,
$$
where $\lambda$ is given by Lemma \ref{step1} (the inequality follows from the fact that $\chi_-\subset\chi$).
With (\ref{couche1}) and (\ref{couche2}) we then get:
$$
\liminf \frac{T_+(B,x)}{\|x\|} \le \liminf \frac{T(B,x)}{\|x\|} \le \limsup \frac{T(B,x)}{\|x\|} 
 \le \limsup \frac{T(B,x+B)}{\|x\|} \le \limsup \frac{T_-(B,x+B)}{\|x\|} 
$$
and then, almost surely:
$$
\lambda(\beta) \le \liminf \frac{T(B,x)}{\|x\|} \le \limsup \frac{T(B,x)}{\|x\|} \le \lambda \le \frac{\lambda(\beta)}{\mu([\beta,+\infty))}.
$$
We finally get:
$$
\lambda \mu([\beta,+\infty)) \le \lambda(\beta) \le \liminf \frac{T(B,x)}{\|x\|} \le \limsup \frac{T(B,x)}{\|x\|} \le \lambda. 
$$
As $\beta$ can be chosen arbitrary small, the lemma is proved. \finpreuve

\bigskip

\proofof{of Theorem \ref{th-vitesse}}
\boulet Because of Lemmas \ref{step1} and \ref{step3}, 
the only thing that remains to be proved about the first part of the theorem is the following convergence in $L^1$:
$$
\lim_{\|x\|\to+\infty} \frac{T(B,x)}{\|x\|}=\lambda.
$$ 
But as the corresponding almost sure convergence holds, as $\frac{T(B,x+B)}{\|x\|}$ converges toward $\lambda$ in $L^1$ and
as we have, for all $x\in\R^d$, $0\le T(B,x) \le T(B,x+B)$, the result follows by uniform integrability arguments.

\boulet Let us prove the second part of the Theorem.
We only consider the case $\lambda>0$ (the other one is easier).
Let $\epsilon \in ]0,\lambda[$.
As $\frac{T(B,x)}{\|x\|}$ converges almost surely toward $\lambda$, there exists an almost surely finite random variable $M$ such that,
for all $x\in\R^d$ such that $\|x\|\ge M$, the following holds:
$$
\lambda-\epsilon\le \frac{T(B,x)}{\|x\|} \le \lambda+\epsilon.
$$ 
Therefore, for all $t\ge0$, we have:
$$
B_{t(\lambda+\epsilon)^{-1}} \setminus B_M \subset S_t \setminus B_M \subset B_{t(\lambda-\epsilon)^{-1}}.
$$
By the second part of Lemma \ref{GD0} we get that almost surely, if $t$ is large enough, then $B_M$ is contained in $S_t$.
Therefore, almost surely, for all large enough $t$ we have:
$$
B_{t(\lambda+\epsilon)^{-1}} \subset S_t \subset B_{t(\lambda-\epsilon)^{-1}}.
$$
This conludes the proof. \finpreuve

\end{document}